\documentclass{article}

\usepackage{latexsym}

\newtheorem{llemma}{Lemma}[section]

\newtheorem{exmp}[llemma]{Example}
\newtheorem{thm}[llemma]{Theorem}

\newtheorem{defn}[llemma]{Definition}

\newtheorem{nnote}[llemma]{Note}

\begin{document}


\title{Similarity measures of intuitionistic fuzzy soft sets and their decision making}

\author{Naim \c{C}a\u{g}man, Irfan Deli \\ \\
          Department of Mathematics, Faculty of Arts and
          Sciences,\\
Gaziosmanpa\c{s}a University, 60250 Tokat, Turkey, \\
             naim.cagman@gop.edu.tr(Naim \c{C}a\u{g}man),\\
             Department of Mathematics, Faculty of Arts and
          Sciences,\\
 Kilis 7 Aral{\i}k University, 79000 Kilis, Turkey, \\
              irfandeli@kilis.edu.tr(Irfan Deli )
       }
\date{}

\maketitle


\begin{abstract}
In this article, we define some types of distances between two
intuitionistic fuzzy soft (IFS) sets and proposed similarity
measures of two IFS-sets. We then construct a decision method which
is applied to a medical diagnosis problem that is based on
similarity measures of IFS-sets. Finally we give two simple example
to show the possibility of using this method for diagnosis of
diseases which could be improved by incorporating clinical results
and other competing diagnosis.

\end{abstract}
\section{Introduction}
In 1999, Molodtsov \cite{mol-99} has introduced the concept of soft
sets. The soft set theory successfully models the problems which
contains uncertainties. In literature, there are theories, such as
probability,  fuzzy sets \cite{zad-65}, intuitionistic fuzzy sets
\cite{ata-86}, rough sets \cite{paw-82} that are dealing with the
uncertain data.

In this work we use soft set theory. The operations (e.g.
\cite{ali-09,cag-09a,maj-03,aslý-11a}) and applications (e.g.
\cite{ali-09,cag-11d,cag-09b}) on soft set theory have been studied
by some researcher. In recent years, many decision making on soft
set theory have been expanded by embedding the ideas of fuzzy sets
(e.g. \cite{ ali-11,ayg-09,cag-09d, deli-12a, deli-12b,deli-12,
fen-08b,fen-10,maj-01a,maju-10,roy-07}), intuitionistic fuzzy sets
(e.g. \cite{ata-86,ata-99,cag-12d, maj-01c, maj-01b,muk-08}) and
rough sets \cite{ ali-11, fen-11}.

Majumdar and Samanta\cite{maju-8} give two types of similarity
measure between soft sets and have shown an application of this
similarity measure of soft sets. Kharal \cite{kha-10} give
counterexamples to show that Definition 2.7 and Lemma 3.5 contain
errors in \cite{maju-8}. In \cite{kha-10}, a new measures have been
presented and this measures have been applied to the problem of
financial diagnosis of firms.

In this paper, we first present the basic definitions and theorem of
soft sets, fuzzy sets, intuitionistic fuzzy sets and intuitionistic
fuzzy soft sets that are useful for subsequent discussions. We then
define distances and similarity measures between two intuitionistic
fuzzy soft (IFS) sets. By using the similarity we construct a
decision making method. We finally give an application, which shows
that the similarity measures can be successfully applied to a
medical diagnosis problem that contains uncertainties.
\section{Preliminary}\label{ss}
In this section, we present the basic definitions of soft set theory
\cite{cag-09a,mol-99}, fuzzy set theory \cite{zad-65},
intuitionistic fuzzy set theory \cite{ata-86} and intuitionistic
fuzzy soft set theory
 \cite{cag-12d} that are useful for subsequent
discussions.
%
%
\begin{defn}\cite{cag-09a}
Let $U$ be a universe, $E$ be a set of parameters that are describe
the elements of $U$, and $A\subseteq E$. Then, a soft set $F_A$ over
$U$ is a set defined by a set valued function $f_A$ representing a
mapping
\begin{equation}\label{soft-set}
f_A: E\to P(U) \textrm{ such that}\, f_A(x)=\emptyset \textrm{ if }
x\in E-A
\end{equation}
where $f_A$ is called approximate function of the soft set $F_A$. In
other words, the soft set is a parametrized family of subsets of the
set $U$, and therefore it can be written a set of ordered pairs
$$
F_A= \{(x, f_A(x)): x\in E, f_A(x)=\emptyset \textrm{ if } x\in
E-A\}
$$
\end{defn}
\begin{defn}\cite{zad-65}
Let $U$ be a universe. Then a fuzzy set $X$ over $U$ is a function
defined as follows:
$$
X=\{(\mu_X(u)/u): u\in U\}
$$
where $\mu_X:U\rightarrow [0.1]$

Here, $\mu_X$ called membership function of $X$, and the value
$\mu_X(u)$ is called the grade of membership of $u\in U $. The value
represents the degree of u belonging to the fuzzy set $X$.
\end{defn}
\begin{defn}\cite{ata-86}
Let $E$ be a universe. An intuitionistic fuzzy set $A$ on $E$ can be
defined as follows:
$$
A =\{<x, \mu_{A}(x),\gamma_{A}(x)>: \,\, x \in E \}
$$
where, $\mu_{A}: E \rightarrow [ 0, 1 ]$ and $\gamma_{A} : E
\rightarrow [ 0, 1 ]$ such that $0 \leq \mu_{A} (x) + \gamma_A(x)
\leq 1$ for any $ x\in E$.

Here, $\mu_{A}(x)$ and $\gamma_{A}(x)$ is the degree of membership
and degree of non-membership of the element $x$, respectively.
\end{defn}

If $A$ and $B$ are two intuitionistic fuzzy sets on $E$, then
\begin{enumerate}
\item $A \subset B$ if and only if
$\mu_{A}(x)\leq \mu_{B}(x)$ and $\gamma_{A}(x) \geq \gamma_{B}(x)$
for $\forall x \in E$

\item $A=B$ if and only if  $\mu_{A}(x)=
\mu_{B}(x)$ and $\gamma_{A}(x) = \gamma_{B}(x)$ $\forall x \in E$

\item $A ^c= \{<x,  \gamma_{A}(x), \mu_{A}(x)>:\,\, x \in E \}$

\item $A\cup B = \{<x, max(\mu_{A}(x),\mu_{B}(x)), min(  \gamma_{A}(x),
\gamma_{B}(x)>:\,\, x \in E \} $,

\item $A\cap B = \{<x, min(\mu_{A}(x),\mu_{B}(x)), max(  \gamma_{A}(x),
\gamma_{B}(x)>:\,\, x \in E \} $,

\item $A + B = \{<x,
\mu_{X}(x)+ \mu_{Y}(x)- \mu_{X}(x)\mu_{Y}(x),
\gamma_{X}(x)\gamma_{Y}(x)>:\,\, x \in E \}$,

\item $A\cdot B
= \{<x, \mu_{A}(x)\mu_{B}(x), \gamma_{A}(x) + \gamma_{B}(x)-
\gamma_{A}(x)\gamma_{B}(x)>:\,\, x \in E \}$.
\end{enumerate}
%
\begin{defn}
  \cite{cag-12d}
 An intuitionistic fuzzy soft set (or namely IFS-set) is
 defined by the set of ordered pairs
$$
\Gamma_A= \{(x, \gamma_A(x)): x\in E, \gamma_A(x)\in \hat{F}(U)\}
$$
where $\gamma_A: E\to \hat{F}(U)$ such that
$\gamma_A(x)=\hat{\emptyset}$ if $x\notin A$ and $\hat{\emptyset}$
is  intuitionistic fuzzy empty  set. Moreover $\gamma_A(x)$ is an
intuitionistic fuzzy set. So it is denoted by
$$\gamma_A(x)=\{(u,\mu_A(u),\nu_A(u)):u\in U\}$$
for all $x\in E$. Moreover, $\mu_A:U\to [0,1]$ and $\nu_A:U\to
[0,1]$ with the condition $0\leq \mu_A(u)+\nu_A(u)\leq 1$, for all
$u\in U$. The numbers $\mu_A(u)$ and $\nu_A(u)$ denote the
membership degree end non-membership degree of $u\in U$ to the
intuitionistic fuzzy set $\gamma_A(x)$, respectively.
\end{defn}
\begin{exmp}\label{ifs-set}
Suppose that there are five car in the universe $U= \{ u_l , u_2,
u_3,u_4,u_5\}$ under consideration  ``$x_1=$large", ``$x_2=$costly",
``$x_3=$secure", ``$x_4=$strong", ``$x_5=$ economic" and", ``$x_6=$
repair". Therefore parameter set is
$E=\{x_1,x_2,x_3,x_4,\\x_5,x_6\}$. Let $A =\{x_1,x_2, x_3, x_4 \}$.
Then IFS-set $\Gamma_A$ is represented the following tabular form;

$$
\begin{array}{rl} \Gamma_A=\bigg\{&
(x_1,\{(u_1,0.5,0.2),(u_2,0.5,0.2),(u_3,0.5,0.2),(u_4,0.5,0.2)\},\\&
(x_2,\{(u_1,0.6,0.4),(u_2,0.9,0.1),(u_3,0.5,0.3),(u_4,0.1,0.9)\},\\&
(x_3,\{(u_1,0.7,0.2),(u_2,0.8,0.1),(u_3,0.2,0.16),(u_4,0.4,0.5)\},\\&
(x_4,\{(u_1,0.4,0.3),(u_2,0.2,0.7),(u_3,0.8,0.2),(u_4,0.2,0.1)\}\bigg\}
\end{array}
$$
\end{exmp}

\section{Similarity Measures of IFS-Sets}

In this section, we first present the basic definitions of distances
between two intuitionistic fuzzy sets \cite{ata-99}  and two soft
sets \cite{maju-8} that are useful for subsequent discussions. We
then define some distances and similarity measures of IFS-sets.

\begin{defn}\cite{ata-99} Let $U = \{x_1, x_2, x_3, . . . , x_n\}$
be a universe and $A, B$ be two intuitionistic fuzzy sets over $U$
with their membership function $\mu_A,\mu_B$ and nonmembership
function $\nu_A,\nu_B$, respectively. Then the distances of $A$ and
$B$ are defined as,

\begin{enumerate}
    \item Hamming distance;
$$
d(A,B)=\frac{1}{2}\sum_{i=1}^n[|\mu_A(x_i)-\mu_B(x_i)|+|\nu_A(x_i)-\nu_B(x_i)|]
$$
    \item Normalized Hamming distance;
$$
l(A,B)=\frac{1}{2n}\sum_{i=1}^n[|\mu_A(x_i)-\mu_B(x_i)|+|\nu_A(x_i)-\nu_B(x_i)|]
$$
    \item Euclidean distance;
$$
e(A,B)=\sqrt{\frac{1}{2}\sum_{i=1}^n[(\mu_A(x_i)-\mu_B(x_i))^2+(\nu_A(x_i)-\nu_B(x_i))^2]}
$$
    \item Normalized Euclidean distance;
$$
q(A,B)=\sqrt{\frac{1}{2n}\sum_{i=1}^n[(\mu_A(x_i)-\mu_B(x_i))^2+(\nu_A(x_i)-\nu_B(x_i))^2]}
$$
\end{enumerate}
\end{defn}

\begin{defn}
\cite{maju-8} Let $U=\{u_1, u_2, u_3,...\}$ be a universe, $E=\{x_1,
x_2, x_3, ...\}$ be a set of parameters, $A,B\subseteq E$, and $F_A$
and $G_B$ be two soft sets on $U$ with their approximate functions
$f_A$ and $g_B$, respectively.

If $A=B$, then similarity between $F_A$ and $G_B$ is defined by
$$
S(F_A, G_B) =\frac{\sum_{i=1}\overrightarrow{f_A(x_i)}\cdot
\overrightarrow{g_B(x_i)}}{\sum_{i=1}
max[\overrightarrow{f_A(x_i)}^2, \overrightarrow{g_B(x_i)}^2]}
$$
where
$$
\overrightarrow{f_A(x_i)}=(
 \chi_{f_A(x_i)}(u_1),\chi_{f_A(x_i)}(u_2),\chi_{f_A(x_i)}(u_3),...)
$$
$$
\overrightarrow{g_B(x_i)}=(
\chi_{g_B(x_i)}(u_1),\chi_{g_B(x_i)}(u_2),\chi_{g_B(x_i)}(u_3),...)
$$
and
$$
\chi_{f_A(x_i)}(u_j)=\left\{\begin{array}{l} {1, u_j\in {f_A(x_i)}}\\
{0, u_j\notin {f_A(x_i)}} \end{array}\right. ,\,\,\,\, \chi_{g_B(x_i)}(u_j)=\left\{\begin{array}{l} {1, u_j\in {g_B(x_i)}}\\
{0, u_j\notin {g_B(x_i)}} \end{array}\right.
$$
\begin{nnote}
If $A\neq B$ and $C =A\cap B\neq\emptyset$, then
$\overrightarrow{f_A(x_i)}=0$ for $x_i\in B/C$ and
$\overrightarrow{g_B(x_i)}=0$ for $x_i\in A/C$.

If $A\cap B=\emptyset$, then $S(F_A, G_B) = 0$ and $S(F_A, F_A^c) =
0$ as $\overrightarrow{f_A(x_i)}\cdot \overrightarrow{f_A^c(x_i)}=0$
for all $i$. \end{nnote}
\end{defn}
\begin{defn}\cite{maju-8}
Let $F_A$ and $G_B$ be two soft sets over $U$. Then, $F_A$ and $G_B$
are said to be $\alpha$-similar, denoted as $F_A \approx^\alpha
G_B$, if and only if $S(F_A, G_B) \geq \alpha$ for $\alpha\in
(0,1).$
\end{defn}
\begin{defn}\cite{maju-8}
Let $U=\{u_1, u_2, u_3, ...\}$ be a universe, $E=\{x_1, x_2, x_3,
...\}$ be a set of parameters, $A,B\subseteq E$ and $F_A,G_B$ be two
soft sets on $U$ with their approximate functions $f_A$ and $g_B$,
respectively. Then, the distances of $F_A$ and $G_B$ are defined as,

\begin{enumerate}
    \item Hamming distance;
$$
d^s(F_A,
G_B)=\frac{1}{m}\bigg\{\sum_{i=1}^m\sum_{j=1}^n|f_A(x_i)(u_j)-g_B(x_i)(u_j)|\bigg\}
$$
    \item Normalized Hamming distance;
$$
l^s(F_A,
G_B)=\frac{1}{mn}\bigg\{\sum_{i=1}^m\sum_{j=1}^n|f_A(x_i)(u_j)-g_B(x_i)(u_j)|\bigg\}
$$
    \item Euclidean distance;
$$
e^s(F_A, G_B)=\sqrt{\frac{1}{m}\sum_{i=1}^m\sum_{j=1}^n(f_A(x_i)(u_j)-g_B(x_i)(u_j))^2}
$$
    \item Normalized Euclidean distance;
$$
q^s(F_A,
G_B)=\sqrt{\frac{1}{mn}\sum_{i=1}^m\sum_{j=1}^n(f_A(x_i)(u_j)-g_B(x_i)(u_j))^2}
$$
\end{enumerate}
\end{defn}

\begin{defn} \cite{maju-8}
Let $F_A$ and $G_B$ be two soft sets over $U$. Then, by using the
Eucledian distance, similarity measure of $F_A$ and $G_B$  is
defined as,
$$
s'(F_A, G_B)=\frac{1}{1+e^s(F_A, G_B)}
$$
Another similarity measure of $F_A$ and $G_B$ can be defined as,
$$
s''(F_A, G_B)=e^{-\alpha e^s(F_A, G_B)}
$$
where $\alpha$ is a positive real number called the steepness
measure.
\end{defn}
\begin{defn}\label{3}
Let $U=\{u_1, u_2, ...,u_n\}$ be a universe, $E=\{x_1, x_2,
...,x_m\}$ be a set of parameters, $A,B\subseteq E$ and
$\Gamma_A,\Lambda_B$ be two IFS-sets on $U$ with their
intuitionistic fuzzy approximate functions
$\gamma_A(x_i)=\{(u,\mu_A(u),\nu_A(u)):u\in U\}$ and
$\lambda_B(x_i)=\{(u,\mu_B(u),\nu_B(u)):u\in U\}$, respectively.

If $A=B$ and $\mu_A(x_i)(u_j)-\nu_A(x_i)(u_j)\not= 0$ or
$\mu_B(x_i)(u_j)-\nu_B(x_i)(u_j)\not= 0$ for at least one $i \in
\{1,2,...,n\}$ and $j \in \{1,2,...,m\}$, then similarity between
$\Gamma_A$ and $ \Lambda_B$ is defined by
\\~\\
$S_{IFS}(\Gamma_A, \Lambda_B)= $
$$
\frac{\sum_{i=1}^m\sum_{j=1}^n|\overrightarrow{(\mu_A(x_i)(u_j)}-\overrightarrow{\nu_A(x_i)(u_j)})\cdot
\overrightarrow{(\mu_B(x_i)(u_j)}-\overrightarrow{\nu_B(x_i)(u_j)})|}{\sum_{i=1}^m\sum_{j=1}^n
max\{\overrightarrow{||\mu_A(x_i)(u_j)}-\overrightarrow{\nu_A(x_i)(u_j)||}^2,
\overrightarrow{||\mu_B(x_i)(u_j)}-\overrightarrow{\nu_B(x_i)(u_j)||}^2\}}
$$
where

$$
\begin{array}{rcl}
\overrightarrow{\mu_A(x_i)(u_j)} & = & (
\mu_A(x_i)(u_1),\mu_A(x_i)(u_2),...,\mu_A(x_i)(u_n))
\\
\overrightarrow{\nu_A(x_i)(u_j)}& = & (
 \nu_A(x_i)(u_1),\nu_A(x_i)(u_2),...,\nu_A(x_i)(u_n))
 \\
 \overrightarrow{\mu_B(x_i)(u_j)}& = &(
 \mu_B(x_i)(u_1),\mu_B(x_i)(u_2),...,\mu_B(x_i)(u_n))
 \\
 \overrightarrow{\nu_B(x_i)(u_j)} & = &(
 \nu_B(x_i)(u_1),\nu_B(x_i)(u_2),...,\nu_B(x_i)(u_n))
\end{array}
$$

If $A=B$ and $\mu_A(x_i)(u_j)-\nu_A(x_i)(u_j)= 0$ and
$\mu_B(x_i)(u_j)-\nu_B(x_i)(u_j)= 0$ for all $i \in \{1,2,...,n\}$
and $j \in \{1,2,...,m\}$, then $S_{IFS}(\Gamma_A, \Lambda_B)=1$.
\end{defn}

\begin{exmp}\label{1} Assume that $U=\{u_1, u_2, u_3,u_4\}$ is a universal set, $E=\{x_1,
x_2,\\
x_3,x_4\}$ is a set of parameters, $A=\{x_1, x_2,x_4\}$, $B=\{x_1,
x_2,x_4\}$ are subsets of $E$. If two IFS-sets $\Gamma_A$ and
$\Lambda_B$ over $U$ are contracted as follows;
$$
\begin{array}{rl}
\Gamma_A= &
\bigg\{(x_1,\{(u_1,0.5,0.5),(u_2,0.4,0.5),(u_3,0.7,0.2),(u_4,0.8,0.1)\}),
\\&(x_2, \{(u_1,0.4,0.6),(u_2,0.2,0.7),(u_3,0.2,0.8),(u_4,0.2,0.2)\}),
\\&(<x_4,\{(u_1,0.2,0.7),(u_2,0.1,0.9),(u_3,0.5,0.4),(u_4,0.7,0.2)\})\bigg\} \\
\end{array}
$$
$$
\begin{array}{rl}
\Lambda_B= &
\bigg\{(u_1,0.2,0.7),(u_2,0.1,0.9),(u_3,0.5,0.4),(u_4,0.4,0.4)\}),
\\&(x_2, \{(u_1,0.5,0.5),(u_2,0.4,0.5),(u_3,0.3,0.6),(u_4,0.4,0.5)\}),
\\&(<x_4,\{(u_1,0.4,0.6),(u_2,0.2,0.7),(u_3,0.2,0.8),(u_4,0.2,0.5)\})\bigg\} \\
\end{array}
$$
Then we can obtain\\

 $\overrightarrow{\mu_A(x_1)(u_j)}=( 0.5,0.4,0.7,0.8) $,\,\, $\overrightarrow{\nu_A(x_1)(u_j)}=(0.5,0.5,0.2,0.1) $,

  $\overrightarrow{\mu_A(x_2)(u_j)}=(0.4,0.2,0.2,0.2) $,\,\, $\overrightarrow{\nu_A(x_2)(u_j)}=(0.6,0.7,0.8,0.2) $,

 $\overrightarrow{\mu_A(x_3)(u_j)}=(0.2,0.1,0.5,0.7) $, \,\,$\overrightarrow{\nu_A(x_3)(u_j)}=(0.7,0.9,0.4,0.2) $,

 $\overrightarrow{\mu_B(x_1)(u_j)}=(0.2,0.1,0.5,0.4) $,\, \,$\overrightarrow{\nu_B(x_1)(u_j)}=(0.7,0.9,0.4,0.4) $,

  $\overrightarrow{\mu_B(x_2)(u_j)}=(0.5,0.4,0.3,0.4) $,\,\, $\overrightarrow{\nu_B(x_2)(u_j)}=(0.5,0.5,0.6,0.5) $,

 $\overrightarrow{\mu_B(x_3)(u_j)}=(0.4,0.2,0.2,0.2) $,\, \,$\overrightarrow{\nu_B(x_3)(u_j)}=(0.6,0.7,0.8,0.5) $.

~\\and\\

$\overrightarrow{(\mu_A(x_1)(u_j)}-\overrightarrow{\nu_A(x_1)(u_j)}=(0.0,-0.1,0.5,0.7)$,

$\overrightarrow{(\mu_A(x_2)(u_j)}-\overrightarrow{\nu_A(x_2)(u_j)}=(-0.2,-0.5,-0.6,0.0)$,

$\overrightarrow{(\mu_A(x_3)(u_j)}-\overrightarrow{\nu_A(x_3)(u_j)}=(-0.5,-0.8,0.1,0.5)$,

$\overrightarrow{(\mu_B(x_1)(u_j)}-\overrightarrow{\nu_B(x_1)(u_j)}=(-0.5,-0.8,0.1,0.0)$,

$\overrightarrow{(\mu_B(x_2)(u_j)}-\overrightarrow{\nu_B(x_2)(u_j)}=(0.0,-0.1,-0.3,-0.1)$,

$\overrightarrow{(\mu_B(x_3)(u_j)}-\overrightarrow{\nu_B(x_3)(u_j)}=(-0.2,-0.5,-0.6,-0.3)$
\\~\\
Now the similarity between $\Gamma_A$ and $\Lambda_B$ is calculated
as
$$
S_{IFS}(\Gamma_A, \Lambda_B)=0.31
$$
\end{exmp}
\begin{thm}
Let $E$ be a parameter set, $A,B\subseteq E$ and $\Gamma_A$ and
$\Lambda_B$ be two IFS-sets over $U$. Then the followings hold;

\begin{enumerate}
    \item[i.] $S_{IFS}(\Gamma_A, \Lambda_B)=S_{IFS}(\Lambda_B, \Gamma_A)$
    \item[ii.] $0\leq S_{IFS}(\Gamma_A, \Lambda_B) \leq 1$
    \item[iii.] $S_{IFS}(\Gamma_A, \Gamma_A)=1$
\end{enumerate}
\end{thm}

\emph{\textbf{Proof:}} Proof easly can be made by using Definition
\ref{3}.

\begin{thm}Let $E$ be a parameter set, $A,B,C\subseteq E$ and $\Gamma_A$, $\Lambda_B$ and
$\Upsilon_C$ be three IFS-sets over $U$ such that $\Gamma_A$ is a
intuitionistic fuzzy  soft subset of $\Lambda_B$ and $\Lambda_B$is a
Intuitionistic fuzzy soft subset of $\Upsilon_C$ then,

$$S_{IFS}(\Gamma_A, \Upsilon_C)\leq S_{IFS}(\Lambda_B, \Upsilon_C) $$

\emph{\textbf{Proof:}}The proof is straightforward.
\end{thm}


\begin{defn}\label{2}
Let $U=\{u_1, u_2, ...,u_n\}$ be a universe, $E=\{x_1, x_2,
...,x_m\}$ be a set of parameters, $A,B\subseteq E$ and
$\Gamma_A,\Lambda_B$ be two IFS-sets on $U$ with their
intuitionistic fuzzy approximate functions
$\gamma_A(x_i)=\{(u,\mu_A(u),\nu_A(u)):u\in U\}$ and
$\lambda_B(x_i)=\{(u,\mu_B(u),\nu_B(u)):u\in U\}$, respectively.
Then the distances of $\Gamma_A$ and $ \Lambda_B$ are defined as,

\begin{enumerate}

\item Hamming distance,
\\~\\
$ d_{IFS}^s(\Gamma_A, \Lambda_B)=$
$$
\frac{1}{2m}\bigg\{\sum_{i=1}^m\sum_{j=1}^n|\mu_A(x_i)(u_j)-\mu_B(x_i)(u_j)|+|\nu_A(x_i)(u_j)-\nu_B(x_i)(u_j)|\bigg\}
$$

\item Normalized Hamming distance,
\\~\\
$l_{IFS}^s(\Gamma_A,\Lambda_B)=$
$$
\frac{1}{2mn}\bigg\{\sum_{i=1}^m\sum_{j=1}^n|\mu_A(x_i)(u_j)-\mu_B(x_i)(u_j)|+|\nu_A(x_i)(u_j)-\nu_B(x_i)(u_j)|\bigg\}
$$

\item Euclidean distance,
\\~\\
$ e_{IFS}^s(\Gamma_A, \Lambda_B)=$
$$
\bigg(\frac{1}{2m}\sum_{i=1}^m\sum_{j=1}^n[(\mu_A(x_i)(u_j)-\mu_B(x_i)(u_j))^2+(\nu_A(x_i)(u_j)-\nu_B(x_i)(u_j))^2]\bigg)^{\frac{1}{2}}
$$


\item Normalized Euclidean distance,
\\~\\
$q_{IFS}^s(\Gamma_A, \Lambda_B)= $
$$
\bigg(\frac{1}{2mn}\sum_{i=1}^m\sum_{j=1}^n[(\mu_A(x_i)(u_j)-\mu_B(x_i)(u_j))^2+(\nu_A(x_i)(u_j)-\nu_B(x_i)(u_j))^2]\bigg)^{\frac{1}{2}}
$$

\end{enumerate}
\end{defn}
%
%
\begin{exmp} \label{6}Let us consider the Example \ref{1}. Then, the
distances of $\Gamma_A$ and $\Lambda_B$ are calculated as follows;
%
%
%
%
%
%
%
$$
\begin{array}{ccc}
d_{IFS}^s(\Gamma_A, \Lambda_B) & = & 0.07 \\
l_{IFS}^s(\Gamma_A, \Lambda_B)& = & 0.37\\
e_{IFS}^s(\Gamma_A, \Lambda_B)& = & 0.28\\
q_{IFS}^s(\Gamma_A, \Lambda_B)& = & 0.19\\
\end{array}
$$

%

\end{exmp}

\begin{thm}Let $E$ be a parameter set, $A,B\subseteq E$ and $\Gamma_A$ and
$\Lambda_B$ be two IFS-sets over $U$. Then the followings hold;
\begin{enumerate}
    \item[i.] $d_{IFS}^s(\Gamma_A, \Lambda_B) \leq n$
    \item[ii.] $l_{IFS}^s(\Gamma_A, \Lambda_B) \leq 1$
    \item[iii.] $e_{IFS}^s(\Gamma_A, \Lambda_B) \leq \sqrt{n}$
     \item[iv.] $q_{IFS}^s(\Gamma_A, \Lambda_B)\leq 1$
\end{enumerate}
\end{thm}

\emph{\textbf{Proof:}} Proof easily can be made by using Definition
\ref{2}.


\begin{thm} Let $IFS(U)$ be a set of all IFS-sets over $U$. Then
the distances functions $d_{IFS}^s$, $l_{IFS}^s$, $e_{IFS}^s$ and
$q_{IFS}^s$, defined from $IFS(U)$ to the non-negative real number
$R^+$, are metric.
\end{thm}

\emph{\textbf{Proof:}} We give only proof for $l_{IFS}^s$. If
$\Gamma_A, \Lambda_B \,and \,\Upsilon_C \in IFS(U)$, then

\begin{itemize}
    \item $l_{IFS}^s(\Gamma_A, \Lambda_B) \leq 0$

$\forall i=\{1,2,...,m\},j=\{1,2,...,n\}$
If
$$
\begin{array}{rl}
  l_{IFS}^s(&\Gamma_A, \Lambda_B)= 0\\& \Rightarrow
|\mu_A(x_i)(u_j)-\mu_B(x_i)(u_j)|+|\nu_A(x_i)(u_j)-\nu_B(x_i)(u_j)|=0
  \\&
\Rightarrow \mu_A(x_i)(u_j)=\mu_B(x_i)(u_j)\wedge
\nu_A(x_i)(u_j)=\nu_B(x_i)(u_j) \\& \Rightarrow
\Gamma_A=\Lambda_B\\
\end{array}
$$
Conversely, let
$$
\begin{array}{rl}
\Gamma_A&=\Lambda_B\\& \Rightarrow
\mu_A(x_i)(u_j)=\mu_B(x_i)(u_j)\wedge
\nu_A(x_i)(u_j)=\nu_B(x_i)(u_j)\\& \Rightarrow
|\mu_A(x_i)(u_j)-\mu_B(x_i)(u_j)|+|\nu_A(x_i)(u_j)-\nu_B(x_i)(u_j)|=0
\\& \Rightarrow l_{IFS}^s(\Gamma_A, \Lambda_B)= 0\\
\end{array}
$$

\item Clearly, $l_{IFS}^s(\Gamma_A, \Lambda_B) = l_{IFS}^s(\Lambda_B, \Gamma_A) $
\item Triangle inequality follows easily from the observation that for any three
IFS-sets $\Gamma_A$. $\Lambda_B$ and $\Upsilon_C$,

$\forall i=\{1,2,...,m\},j=\{1,2,...,n\}$

$|\mu_A(x_i)(u_j)-\mu_B(x_i)(u_j)|+|\nu_A(x_i)(u_j)-\nu_B(x_i)(u_j)|=|\mu_A(x_i)(u_j)-\mu_C(x_i)(u_j)+\mu_C(x_i)(u_j)-\mu_B(x_i)(u_j)|+
|\nu_A(x_i)(u_j)-\nu_C(x_i)(u_j)+\nu_C(x_i)(u_j)-\nu_B(x_i)(u_j)|\leq
|\mu_A(x_i)(u_j)-\mu_C(x_i)(u_j)|+|\mu_C(x_i)(u_j)-\mu_B(x_i)(u_j)|+|\nu_A(x_i)(u_j)-\nu_C(x_i)(u_j)|+|\nu_C(x_i)(u_j)-\nu_B(x_i)(u_j)|$

Therefore, we have:

 $l_{IFS}^s(\Gamma_A, \Lambda_B)
\leq l_{IFS}^s(\Gamma_A, \Upsilon_C) + l_{IFS}^s(\Upsilon_C,
\Lambda_B) $
\end{itemize}

The others proofs can made similarly.

\begin{defn}\label{4}
Let $\Gamma_A$ and $\Lambda_B$ be two IFS-sets over $U$. Then, by
using the Hamming distance, similarity measure of $\Gamma_A$ and
$\Lambda_B$ is defined as,

$$
S'_{IFS}(\Gamma_A, \Lambda_B)=\frac{1}{1+d_{IFS}^s(\Gamma_A,
\Lambda_B)}
$$

Another similarity measure of $F_A$ and $G_B$ can be defined as,
$$
S''_{IFS}(\Gamma_A, \Lambda_B)=e^{-\alpha d_{IFS}^s(\Gamma_A,
\Lambda_B)}
$$
where $\alpha$ is a positive real number called the steepness
measure.
\end{defn}
\begin{defn}
Let $\Gamma_A$ and $\Lambda_B$ be two IFS-sets over $U$. Then,
$\Gamma_A$ and $\Lambda_B$ are said to be $\alpha$-similar, denoted
as $\Gamma_A \approx^\alpha \Lambda_B$, if and only if $S'(\Gamma_A,
\Lambda_B) \geq \alpha$ for $\alpha\in (0,1).$

We call the two IFS-sets significantly similar if
$S'_{IFS}(\Gamma_A, \Lambda_B)
> \frac{1}{2}$.
\end{defn}

\begin{exmp}Let us consider the Example \ref{6}.
Similarity measure of $\Gamma_A$ and $\Lambda_B$ is obtained as,

$$
S'_{IFS}(\Gamma_A, \Lambda_B)=\frac{1}{1+d_{IFS}^s(\Gamma_A,
\Lambda_B)}=0.73
$$

$\Gamma_A$ and $\Lambda_B$ is significantly similar because
$S'_{IFS}(\Gamma_A, \Lambda_B)=0.73> \frac{1}{2}$
\end{exmp}
\begin{thm} Let $E$ be a parameter set, $A,B\subseteq E$ and $\Gamma_A$ and
$\Lambda_B$ be two IFS-sets over $U$. Then the followings hold;
\begin{enumerate}
    \item[i.] $0\leq S'_{IFS}(\Gamma_A, \Lambda_B) \leq 1$
    \item[ii.]  $S'_{IFS}(\Gamma_A, \Lambda_B)=S'_{IFS}(\Lambda_B, \Gamma_A)$
    \item[iii.]  $S'_{IFS}(\Gamma_A, \Lambda_B)=1 \Leftrightarrow \Gamma_A= \Lambda_B$
\end{enumerate}
\end{thm}
\emph{\textbf{Proof:}} Proof easly can be made by using Definition
\ref{4}.\\
%
\section{Decision Making Method}
In this section, we construct a decision making method that is based
on the similarity measure of two IFS-sets. The algorithm of decision
making method can be given as;
\begin{description}
    \item[\emph{Step 1.}] Constructs a IFS-set $\Gamma_A$ over $U$ based on an expert,
    \item[\emph{Step 2.}] Constructs a IFS-set $\Lambda_B$ over $U$ based on a responsible person for the problem,
    \item[\emph{Step 3.}] Calculate the distances of $\Gamma_A$ and $\Lambda_B$,
    \item[\emph{Step 4.}] Calculate the similarity measure of $\Gamma_A$ and $\Lambda_B$,
    \item[\emph{Step 5.}] Estimate result by using the similarity.
\end{description}

~\\

Now, we can give an application for the decision making method. By
using the Hamming distance, similarity measure of two IFS-sets can
be applied to detect whether an ill person is suffering from a
certain disease or not.

\section{Application}
In this applications, we will try to estimate the possibility that
an ill person having certain visible symptoms is suffering from
cancer. For this, we first construct a IFS-set for the illness and a
IFS-set for the ill person. We then find the similarity measure of
these two IFS-sets. If they are significantly similar, then we
conclude that the person is possibly suffering from cancer.
\begin{exmp}\label{exm-1}
Assume that our universal set contain only two elements cancer and
not cancer, i.e. $U=\{u_1, u_2\}$. Here the set of parameters
$A=B=E$ is the set of certain visible symptoms, let us say,
$E=\{x_1, x_2, x_3,x_4,x_5, x_6, x_7,x_8,x_9 \}$ where $x_1=$
jaundice, $x_2=$ bone pain, $x_3=$ headache, $x_4=$ loss of
appetite, $x_5=$ weight loss, $x_6=$ heal wounds , $x_7=$ handle and
shoulder pain, $x_8=$ lump anywhere on the body for no reason and
$x_9=$ chest pain.
\end{exmp}

\begin{description}
\item[\emph{Step 1.}] Constructs a IFS-set $\Gamma_A$ over $U$ for cancer is given below and this can be
prepared with the help of a medical person:

$$
\begin{array}{rl}
\Gamma_A=
\big\{&(x_1,\{(u_1,0.5,0.5),(u_2,0.4,0.5)\},(x_2,\{(u_1,0.7,0.2),(u_2,0.8,0.1)\}),
\\&(x_3, \{(u_1,0.4,0.6),(u_2,0.2,0.7)\},(x_4,\{(u_1,0.2,0.8),(u_2,0.2,0.2)\}),
\\&(x_5,\{(u_1,0.2,0.7),(u_2,0.1,0.9)\},(x_6,\{(u_3,0.5,0.4),(u_4,0.7,0.2)\}),
\\&(x_7,\{(u_1,0.3,0.7),(u_2,0.4,0.4)\},(x_8,\{(u_1,0.5,0.2),(u_2,0.7,0.1)\}),
\\&(x_9,\{(u_1,0.3,0.4),(u_2,0.7,0.1)\})\,\,\,\,\,\big\} \\
\end{array}
$$

\item[\emph{Step 2.} ] Constructs a IFS-set $\Lambda_B$ over $U$ based on data of ill
person:
$$
\begin{array}{rl}
\Lambda_B=
\bigg\{&(x_1,\{(u_1,0.9,0.1),(u_2,0.9,0.0)\},(x_2,\{(u_1,0.1,0.9),(u_2,0.1,0.8)\}),
\\&(x_3, \{(u_1,0.7,0.1),(u_2,0.8,0.9)\},(x_4,\{(u_1,0.9,0.1),(u_2,0.9,0.8)\}),
\\&(x_5,\{(u_1,0.9,0.1),(u_2,0.9,0.2)\},(x_6,\{(u_3,0.1,0.9),(u_4,0.1,0.8)\})
\\&(x_7,\{(u_1,0.9,0.1),(u_2,0.7,0.9)\},(x_8,\{(u_1,0.9,0.9),(u_2,0.1,0.9)\}),
\\&(x_9,\{(u_1,0.8,0.1),(u_2,0,1)\})\,\bigg\} \\
\end{array}
$$

\item[\emph{Step 3.} ] Calculate Hamming distances of $\Gamma_A$ and $\Lambda_B$,

$$d_{IFS}^s(\Gamma_A, \Lambda_B)\cong 1.1$$
\item[\emph{Step 4.} ] Calculate the similarity measure of $\Gamma_A$ and $\Lambda_B$,

$$
 S'_{IFS}(\Gamma_A, \Lambda_B)=\frac{1}{1+d_{IFS}^s(\Gamma_A,
\Lambda_B)}\cong 0.48< \frac{1}{2 }
$$
%
\item[\emph{Step 5.} ] Hence the two IFS-sets, i.e. two symptoms $\Gamma_A$ and $\Lambda_B$
are not significantly similar. Therefore, we conclude that the
person is not possibly suffering from cancer.
\end{description}
\begin{exmp}
Let us consider Example \ref{exm-1} with different ill person.
\end{exmp}
\begin{description}
\item[\emph{Step 1.} ] Constructs a IFS-set for cancer $\Gamma_A$ is in the Example \ref{exm-1}:

\item[\emph{Step 2.} ] A person suffering from the following symptoms whose corresponding
IFS-set $\Upsilon_C$ is given below:
$$
\begin{array}{rl}
\Upsilon_C=
\bigg\{&(x_1,\{(u_1,0.5,0.4),(u_2,0.4,0.4)\},(x_2,\{(u_1,0.7,0.1),(u_2,0.8,0.1)\}),
\\&(x_3, \{(u_1,0.4,0.5),(u_2,0.2,0.6)\},(x_4,\{(u_1,0.2,0.7),(u_2,0.2,0.1)\}),
\\&(x_5,\{(u_1,0.2,0.6),(u_2,0.1,0.8)\},(x_6,\{(u_3,0.5,0.3),(u_C,0.7,0.1)\})
\\&(x_7,\{(u_1,0.2,0.6),(u_2,0.1,0.8)\},(x_8,\{(u_1,0.5,0.3),(u_2,0.7,0.1)\}),
\\&(x_9,\{(u_1,0.5,0.3),(u_2,0.7,0.1)\})\bigg\} \\
\end{array}
$$

\item[\emph{Step 3.} ] Calculate Hamming distances of $\Gamma_A$ and $\Lambda_B$,
$$d_{IFS}^s(\Gamma_A, \Lambda_B)\cong 0,41$$
 \item[\emph{Step 4.} ] Find the similarity measure of these two IFS-sets as:
$$
S'_{IFS}(\Gamma_A, \Upsilon_C)=\frac{1}{1+d_{IFS}^s(\Gamma_A,
\Upsilon_C)}\cong 0.71>\frac{1}{2}
$$

 \item[\emph{Step 5.} ] Here the two IFS-sets, i.e. two symptoms $\Gamma_A$ and $\Upsilon_C$
are significantly similar. Therefore, we conclude that the person is
possibly suffering from cancer.
\end{description}

\section{Conclusion}
Majumdar and Samanta\cite{maju-8} give two types of similarity
measure between soft sets and have shown an application of this
similarity measure of soft sets. In \cite{kha-10}, Kharal give
counterexamples to show that Definition 2.7 and Lemma 3.5 contain
errors in \cite{maju-8}. In \cite{kha-10}, a new measures have been
presented and this measures have been applied to the problem of
financial diagnosis of firms. In this paper, we have defined four
types of distances between two IFS-sets and proposed similarity
measures of two IFS-sets. Then, we construct a decision making
method based on the similarity measures. Finally, we give two simple
examples to show the possibility of using this method by using
Hamming distance for diagnosis of diseases. In these example, if we
use the other distances, we can obtain similar result.

The method can be applied to problems that contain uncertainty such
as problems in social, economic systems, pattern recognition,
medical diagnosis, game theory, coding theory and so on.


\end{document}